\documentstyle[12pt,twoside,leqno]{article}
\title {Defect and evaluations}
\author{Stefano Trapani}
\date  {} 

\newtheorem{theorem}             {Theorem}    [section]
\newtheorem{corollary}  [theorem]{Corollary}
\newtheorem{lemma}      [theorem]{Lemma}
\newtheorem{proposition}[theorem]{Proposition}
\newtheorem{definition}    [theorem]{Definition}
\newtheorem{examples}    [theorem]{Examples} 
\newtheorem{remark}    [theorem]{Remark}      %
\newenvironment{proof}{\begin{trivlist}\item[]{\em Proof.}}%
                      {\rule{.1in}{.1in}\end{trivlist}}
\newcounter    {assertcount}

\newcommand{\C} {{\bf  C}}

\newcommand{\R} {{\bf  R}}
\newcommand{\Z} {{\bf  Z}}

\newcommand{\de}   {\partial}

\newcommand{\Calpha} {C^{\alpha}( \partial \Delta, {\bf R}^N)} 
\newcommand{\LA} {{ \cal L}_A}

\newcommand{\EA} {{\cal E}_A}
\newcommand{\EAC} {{{\cal E}_A}^{\C}}
\newcommand{\ONalpha}   {{{\cal O}^N}_{\alpha} }

\renewcommand{\epsilon}{\varepsilon}

\begin{document}

\maketitle

\begin{abstract}
Let $S$ be a generic submanifold of $\C^N$ of real codimension m. 
In this work we continue  the study, carried over by various authors, of the 
set of analytic discs  attached to S. Moreover we look at the subspaces of 
$\C^N$ obtained by evaluating at $0$ and at $1,$ holomorphic  maps which are 
the infinitesimal deformations of analytic discs attached to $S.$ \end{abstract}

\section{Introduction} \label{Int}

Let $S$ be a real submanifold of $\C^N$ of real codimension $m,$
assume that $S$  is generic i.e. that  $T_p(S)+i T_p(S)= \C^N,$
for every $p \in S.$
($T_p(S)$ is the tangent space to $S$ at $p$). In the study of the 
polynomial hulls of such submanifolds, as well as in the study of 
holomorphic extendability of their C-R functions,  it is of particular 
importance the set of analytic discs attached to the submanifold.
Let $\Delta$ be the unit disc $ |\zeta| < 1 $  in $\C.$
A continuous map $f: \bar{\Delta} \rightarrow \C^N$ is said to be an 
analytic disc attached to $S,$ if $f(\de \Delta) \subseteq S$ and $f$ is 
holomorphic on $\Delta.$
Let $\alpha$ be a positive real which is not an integer, 
it is convenient to consider analytic discs attached to $S$ which are in 
$C^{\alpha}(\bar{\Delta},\C^N).$ 
If $S$ is given by the zero set of a maximal rank defining map 
$\rho : \Omega  \rightarrow \R^m,$ ($\Omega$ being an open neighborhood of 
$S$ in $\C^N$), and $f_0 : \bar{\Delta} \rightarrow \C^N $ is an analytic 
disc attached to $S,$ then the set $M$ of such analytic discs is  
the zero set of the map $ f \rightarrow \rho \circ f.$
Let  $\ONalpha$ be the space of $C^{\alpha}$ maps from  $\de \Delta$ to 
$\C^N$ which extends holomorphically to $\bar{\Delta}.$
The defining map of $M$ goes from  $\ONalpha$  to $\Calpha,$ and its 
differential is  given by 
\[ w \rightarrow Re \left(\frac{\de \rho}{\de z} \circ f_0 (w(\zeta))
\right) \]
Because of the implicit function theorem, the set $M$ is a manifold, 
as soon as its differential is onto.   
Whenever $M$ is a manifold, its tangent space at $f_0$ is given by
\[ \{ w \in \ONalpha : Re \left(\frac{\de \rho}{\de z}(f_0(\zeta))w
\right)=0 \}.\]
Same procedure can be carried over in the case we want to look at the set of 
analytic discs with a fixed point at the boundary. In other words the set  
$ M_q= \{ f \in M : f(1)=q \},$ where $q=f_0(1).$
In \cite{Ce},  \cite{Fo}, \cite{Gl}, \cite{Oh}, and \cite{Trep}, sufficient 
condition are given for either $M$ or $M_q$ to be manifolds in a neighborhood 
of $f_0.$  
\vskip .1cm
To gather informations on the  polynomial 
hull of $S,$ it is obviously interesting to describe the union of the images of 
$f(\Delta)$ for $f$ in $M $ and $M_q.$ Or, which is the same because of 
the  action on $M$ of the group of authomorphims of $\Delta,$ the images of 
the evaluation map $\mu_{0} : M \rightarrow \C^N $ given by $\mu_{0}(f)=f(0).$ 
To derive informations on  
the evaluation maps, one  might study  the image of their differentials 
(i.e. the maps themselves), restricted to the tangent space to $M $ and 
$M_q.$ This will tell us when to apply the open mapping theorem. 
\vskip .1cm
Similarly, the evaluation maps at the boundary,  
$\mu_{1}(f)=f(1)$ on $M,$ and  $\mu_{-1}(f)=f(-1)$ on $M_q,$ are 
related to the local extendability  of $C-R$ functions on $S,$ 
see \cite{Tu}. 
\vskip .1cm
In this paper we will give a sufficient criteria for $M_q$ to be a 
manifold (more general then the one stated in \cite{Ce}). We will also 
describe the images of the differentials of the above evaluation 
maps, restricted to the suitable tangent spaces. 
We will 
determine their dimensions. (Some work in this direction was previously made 
in \cite{Ce}, \cite{Eisen}, \cite{Pa}, and \cite{Pa2}). 
In particular we show that, for $x \in \Delta,$  the dimension of the image of
the differential of $ {\mu_x}_{|M}$ and  $ {\mu_x}_{|M_q}$ are independent 
on $x.$ The same is true as $x$ varies in $\de \Delta.$ 
More, we prove that, as $x$ varies in $\bar{\Delta}-1,$ the complex 
span of the image of the differential of ${\mu_{x}}_{|M_q}$ 
defines a $C^{\alpha}$ vector  bundle which is 
holomorphic on $\Delta.$ Moreover its restriction to $\de \Delta$ contains 
the complex tangent bundle to $S$ as a subbundle.   
We will then apply all this, as outlined above, 
to find the images of the evaluation maps themselves.   
For example we show that, if every point of S is minimal, then there is an 
open dense subset $\Omega$ of $M,$ invariant under the action of the 
authomorphism group of $\Delta,$  such that the restriction to $\Omega$  
of the evaluation $\mu_0$ is an open map. 
In the present work, we will make heavy use of the results of \cite{Trep} 
about the operator 
\[ U_A(w)= Re \left(\frac{\de \rho}{\de z}(f_0(\zeta))w  \right) \] 
from $\ONalpha$ to $\Calpha,$ whose kernel is the tangent space to $M$ at 
$f_0.$  
What was proved in \cite{Trep},  was that the study of the above operator 
can be reduced to the study of its continuous $L^2$ extension which is much 
easier to deal with. We will mimic this procedure and apply it to the 
restriction of $U_A$ to suitable subspaces of $\ONalpha.$ 
The Theorems about the evaluation maps will follow fairly easily.

\vskip .1cm
In the last two sections we consider the special case, interesting and 
substantially simpler, of totally real submanifolds, and of hypersurfaces.
\vskip .1cm
The Author would like to thank Prof. Claudio Rea and Doctor Patrizia Rossi, 
for introducing him to the subject, and for many stimulating discussions.

\section{Preliminaries} \label{prel}

Let $\alpha$ be a positive real which is not an integer, write 
$\alpha= k + \beta,$ where $k$ is a non negative integer and $0 < \beta < 1.$  
Denote by $C^{\alpha}$ the space of maps which are continuously differentiable 
up to order k, and such that all their derivatives of order  k, are Holder 
maps with Holder coefficient $\beta.$ We assume $\beta > 1/2.$
Let  $\Delta$ be the disc  $ |\zeta| < 1$ in $\C.$ 
Let  $\ONalpha$ be the space of $C^{\alpha}$ maps from  $\de \Delta$ to 
$\C^N$ which extends holomorphically to $\bar{\Delta}.$
Denote with Re z and  Im z respectively the real and the imaginary part 
of a complex number $z.$
Let $f_0$ be an analytic disc attached to a generic $C^2$ submanifold  $S$ of 
real codimension $m.$
Then $f_0^*(T(S))$ is a vector bundle over $\de \Delta.$
If $f_0^*(T(S))$ is trivial over $\de \Delta,$ (orientable case),  
there exists an open 
neighborhood U of the graph of $f_0$ in $\de \Delta \times \C^N,$ and a 
$C^2$ map $\rho: U \rightarrow \R^m$ such that for every fixed $x \in \de 
\Delta$ the map $Y \rightarrow \rho(x,Y)$ is a defining map of $S.$ 
If $f_0^*(T(S))$ is non trivial, then $(f_0 \circ sq)^*(T(S))$ is.
Here $sq(\zeta)= \zeta^2,$ see \cite{Trep}. 
All we say in this paper can be easily carried over to the non 
orientable case,  by composing with the map $sq$ and modifying the function 
spaces accordingly, as in \cite{Trep}. However  we will for simplicity confine 
ourself to the orientable case. So for us every analytic disc attached to 
S is an orientable one.
Observe that the map from $\ONalpha$ to $C^{\alpha}(\de \Delta,\R^m)$ given by 
$ f \rightarrow \rho \circ f$ is $C^1.$ (see appendix in \cite{Gl}.) 
Since S is generic, then the matrix   
$ \frac{\de \rho}{\de z}(\zeta,f_0(\zeta))$ has maximal rank $m$ for 
every $\zeta \in \de \Delta.$
We will consider this problems sketched in the Introduction from a slightly 
more general point of view:   

Let  $A : \de \Delta \rightarrow {\cal M}(m \times N, \C)$ a 
$C^{\alpha}$ map,  where   ${\cal M}(m \times N, \C)$ denotes the space of 
$ m \times N $ complex matrices.  Assume that  $ m \leq N $ and that $A(\zeta)$ 
has maximal rank $m$ for every $ \zeta \in \de \Delta.$
(In the case of an analytic disc $f_0$ the map $A$ is given by   
$A(\zeta)= \frac{\de \rho}{\de z}(\zeta,f(\zeta)).$
Let 
\[ { \cal L}_A= \{ w \in { \cal O}^N : Re(A w)=0 \ \mbox{\ on \ } \de \Delta 
\}. \]
and let   
\[ \EA= \{ \gamma  \in C^{\alpha}( \de \Delta, \R^N) \] 
\[ \mbox{such that} \  \gamma^t A \ \mbox{extends holomorphically to} 
 \ \overline{\Delta} \}. \]
Let   $ \EAC \ \mbox{be the complexification of} \ \EA,$
that is 
\[ \EAC= \{ \delta  \in C^{\alpha}( \de \Delta, \C^N) : Re \delta 
\ \mbox{and  \ }  Im \delta \ \mbox{ \ belong to \ } \EA \}. \]
For a given $ x \in \overline{\Delta}. $ Define
the evaluation map
\[{\mu}_x :{\cal O}^N  \rightarrow \C^N  \]
as the value in $x$ of the holomorphic extension of $w.$ 
Define
\[ {\Psi}_x : \EAC \rightarrow \C^N  \ 
\mbox{as the evaluation at $x$ of the holomorphic extension of $\delta^t A.$
}\] In other words 
\[\begin{array}{l} {\Psi}_x(\delta)= \frac{1}{2 \pi i} \int_{\de \Delta} 
\frac{\delta^t A(\zeta)}{\zeta - x} d \zeta \\
\mbox{for $x \in \Delta$ and} \\ 
\Psi_x(\delta)= \delta^t(x) A(x) \\
\mbox{for $ x \in \de \Delta.$} \end{array} \]
By the results in \cite{Pa2} the space   $\EA$ has a finite 
dimension $d.$ The number $d$ will be called defect of  $A.$ 
Let  $ \gamma_1 \ldots \gamma_{d} $ be a  base of    $\EA$ over $\R,$
so  $ \gamma_1 \ldots \gamma_{d} $ is also a base of  
$ \EAC \ \mbox{ \ over \ } \C.$
Given $x \in \overline{\Delta},$ if  $\EA \neq 0 $ define $P(x)$ as  the 
matrix  
having rows ${\Psi}_x(\gamma_1), \ldots {\Psi}_x(\gamma_{d}).$ If 
$\EA=0$ let  $P(x)$ be the zero matrix in  ${\cal M}(N \times N, \C).$   

Given $x \in \bar{\Delta,}$
denote by $(K_A)_x$ the  kernel of the map ${\Psi}_x$ and by  $(N_A)_x$ 
the kernel of its restriction to $\EA.$ We will be especially interested in 
$(K_A)_0$ and $(N_A)_0$ which will be simply  denoted by $K_A$ and $N_A$
respectively.  
Let $l_x$  be the real dimension of the image of  
${{\Psi}_x}_{|{\cal E}_A} $, and let  $r_x$ be the complex 
dimension of the image of ${\Psi}_x$, i.e. the rank of $P(x).$

Finally let $Q_A$ be  the subspace of $\EAC$  given by
\[ Q_A = \{ \delta \in \EAC : \bar{\delta} = - \bar{\zeta} \delta. \} \]
Note that if $\delta \in Q_A,$ then $\delta^t A = - \zeta \bar{\delta}^t 
A,$ therefore the holomorphic extension of $\delta^t A$ vanishes at $0,$ 
i.e. $Q_A \subseteq K_A.$

We can relate $Q_A$ and $\EA,$ via the following 

\begin{lemma}

Let  $\gamma$ be in  $\EA$ such that   $\gamma(1)=0,$ denote by  $\delta$  
the map 
\[ \delta= \frac{\gamma}{1-\bar{\zeta}}. \] Then $\delta$ 
is in $C^{\alpha}(\de \Delta,\C^m), $ more precisely $\delta \in Q_A.$
Viceversa, if $\delta$ is in $Q_A,$ then $(1- \bar{\zeta}) \delta$ is in 
$\EA$ and it vanishes at $1.$ \label{division} \end{lemma}

\begin{proof}

Given $ \delta \in Q_A, $ then $(1- \bar{\zeta}) \delta$ is real and 
vanishes at $1.$ Moreover $(1- \bar{\zeta}) \delta^t A= \delta^t A - 
\frac{1}{\zeta} \delta^t A $ which extends holomorphically to 
$\bar{\Delta}$ since  $Q_A \subseteq  K_A \subseteq \EAC.$   
Viceversa
Let  $\gamma \in \EA$  such that  $\gamma(1)=0.$  
Since $ \beta > 1/2 $ we can write $A=(A',A")$ where $A'$ has 
values in $GL(m,\C),$ see \cite{Trep}. We have the factorization
$- {A'}^{-1} \bar{A'}= \Theta^{-1} \Lambda \bar{\Theta}$   
 where  $ \Theta $ is a holomorphic map from  $\bar{\Delta}$ to 
$GL(m,\C),$ and $\Lambda$ is a diagonal matrix of the form 
\[ \Lambda = \left( \begin{array}{llcr} \zeta^{k_1} & 0 &  \ldots & 0 \\   
0 & \zeta^{k_2} & \ldots & 0 \\ \ldots & \ & \ & \   \\ 
0 & \dots & \ &  \zeta^{k_m} \\ \end{array} \right). \] 
Here $k_1 \ldots k_m$ are suitable integers, called partial indices of  
$A,$  (Birkhoff factorization), See \cite{Gl}, \cite{Pa} and \cite{Ve}.
Now $\gamma^t A'$ extends holomorphically to $\bar{\Delta}.$
Set $\gamma^t A'= u^t,$ then  $ u^t {A'}^{-1} $ is real, i.e.
$ - u^t \Theta^{-1} \Lambda \bar{\Theta}= \bar{u^t}.$
By setting  $v^t = i u^t \Theta^{-1},$ we find   $v_j= \bar{v_j} \zeta^{-k_j}$ 
on $\de \Delta,$ for $1 \leq j \leq m.$  In particular each component  $v_j$ 
is a polynomial vanishing at $1.$ Therefore we have $v_j=(1- \zeta) w_j $ 
with  $w_j$ polynomial. We then obtain   
\[\bar{\delta}= \frac{\gamma}{1-\zeta}= -i ({{A'}^t})^{-1} \Theta^{-1} w. \] 
Where $w$ is the  vector with components $w_j.$ Hence   $\delta$ and  
$\bar{\delta}$   are in   $C^{\alpha}(\de \Delta,\C^m).$ 
Let $f$ be the holomorphic extension of 
$ \gamma^t A $ to $\bar{\Delta}.$ Let 
$f_n = \frac{f}{1-\zeta + 1/n}.$ The maps  
$f_n$ are holomorphic in $\bar{\Delta}$ and converge 
uniformly to  $ \frac{f}{1-\zeta} = \bar{\delta}^t A$ on $\de \Delta.$
By the maximum principle the sequence $f_n$ converges  uniformly to a 
holomorphic map on $\bar{\Delta}$ which is an extension of $\bar{\delta}^t A.$
Since $\gamma$ is real, $ \delta = -\zeta \bar{\delta}.$
Hence $\delta^t A $ extends holomorphically to $\bar{\Delta}, $ therefore 
$ \delta \in Q_A.$ \end{proof}

\begin{lemma}

The space  $Q_A$ is maximal totally real in $K_A.$ \label{totreal} 
\end{lemma}

\begin{proof}

Clearly  $Q_A \cap i Q_A =0.$ Given $\delta \in K_A$ we have  
$\delta = \delta' +i \delta" $ with 
$\delta'= \frac{ \delta + \zeta \bar{\delta}}{2},$ 
and $\delta"= \frac{ \delta - \zeta \bar{\delta}}{2i}.$ \end{proof}

\begin{remark}
There exist two   
natural embeddings of $N_A$ into $Q_A.$ 
One is given by $ \gamma \rightarrow  (1- \zeta) \gamma $ and its image 
coincides with the set of elements in $Q_A$ vanishing at $1.$ 
The other is given by  $ \gamma \rightarrow  i (1+ \zeta) \gamma $
and its image coincides with the set of elements in $Q_A$ vanishing at $-1.$ 
The proof is as in Lemma \ref{division}. \label{division2} \end{remark} 

\begin{proposition}

\[ \begin{array}{l} a) \  \mbox{The numbers  $r_x$ are independent} \\ 
\mbox{on $x$ for  $x$ in $\bar{\Delta}.$}
\ \mbox{Let $r$ be their common value.}
\\
b) \ \mbox{The numbers  $l_x$ are independent on $x$ for  $x$ in 
$\Delta.$ Let $l$ be their common value.} \\
c) \ \mbox{We have $l_x= r$ for every  $x \in \de \Delta,$} 
\end{array} \] \label{ranks} \end{proposition}

\begin{proof}
Proof of a) 
We will show that $r_1=r_0.$  
Given any $x \in \Delta,$ let $\sigma_x$ the authomorphism of the disc 
sending $0$ to $x$ and keeping $1$ fixed, it will be sufficient to apply 
the above equality to $A \circ \sigma_x$ to show that $r_x$ is 
independent on $x$ for $x$ in $\Delta.$ By a further 
application to $A \circ \sigma,$ where $\sigma$ is an rotation, we conclude 
that $r_x$ is independent on $x \in \bar{\Delta.}$ 
Since $A(1)$ has maximal rank m, then an element $\gamma$ in  
$\EAC$ is in the kernel of $\Psi_{1}$ if and only if  
$\gamma= \gamma_1 +i \gamma_2$ with  $\gamma_1$ and  $ \gamma_2$ in  
$\EA$ vanishing at  $1.$ It follows from Lemma  \ref{division},  
and Lemma \ref{totreal}, that the  kernel of $\Psi_{1}$ has real dimension  
$d- r_0,$ i.e.  $r_1 =r_0.$
So we proved part a). Part c) follows directly from the fact that $A(x)$ 
has rank $m$ for $x \in \de \Delta.$

Let us prove b)
Given $x$ in $\Delta$ different from $0,$ let  
\[ t(\zeta)=\frac{\bar{x}{\zeta}^2 -(|x|^2+1)\zeta +x}{\zeta}, \] the function  
$t(\zeta)$ has only one pole in $\bar{\Delta}$ and it is a simple pole at  
$0.$ It also has only one zero in $\bar{\Delta}$ and it is a simple zero at $x.$
Moreover  $t(\zeta)$ is real on $\de \Delta.$ 
Hence the map $ \gamma \rightarrow t(\zeta) \gamma(\zeta)$ defines a linear 
isomorphism  from $(N_A)_0$ onto $(N_A)_x.$ \end{proof}  

\begin{remark}

If  $\gamma$ is  in ${{\cal E}_A}$ and  $w$ is  in $ { \cal L}_A$ then the 
map  
\[ x \rightarrow {\Psi}_x(\gamma) \mu_x(w) \] 
is holomorphic in $\Delta$   and it is purely imaginary on $\de \Delta.$
Hence it is a purely imaginary constant.
\label{pearing} \end{remark}

\begin{proposition}

The integers  $r$ and  $l$ have the following properties 
\[ \begin{array}{l} a)  \  0 \leq r \leq min(l,m), \  l \leq min(d,2r),
  \ r=0 \ \mbox{if and only if} \ d=0  \\ 
b) \  \mbox{ \ If \ } r=m  \ \mbox{then \ } \\  
\LA= \{ w \in \ONalpha : \mbox{such that $P(x)w(x)$ is a purely imaginary
constant}  \}. \\
c)  \  \mbox{If the map  A is associated to a small analytic disc}  \\
\mbox{attached to a generic submanifold S, then  $r=d.$} \end{array} \] 
\label{properties}  \end{proposition}

\begin{proof}

The inequalities $ r \leq l \leq 2 r,$ and $l \leq d,$  
follows immediately from the definitions.   
If $r=0$ then $d=0$ because of Proposition \ref{ranks}.
Fix a point $ x \in \de \Delta, $ and a base 
$\gamma_1, \ldots  \gamma_{d}$ of $\EA.$  
Let $M(x)$ be the matrix having the vectors 
$\gamma_1(x), \ldots \gamma_{d} (x),$ as rows,
so $P(x)= M(x)A(x).$  
Since $A(x)$ has rank m,  the matrix  $M(x)$ must have rank  $r,$  
hence $r \leq m.$  

Assume now that $r=m.$  Given $w \in \LA,$ we know from Remark 
\ref{pearing} that $P w$ is a purely imaginary constant. 

Fix viceversa a map  $ w \in \ONalpha$ such that $Pw$ is a purely imaginary 
constant.    
We have that $M(x)(Re A(x)w(x))=0,$ for $x \in \de \Delta.$ 
Since $ r=m,$ then $M(x)$ defines a 1-1 linear map from $\R^{d}$ into 
$\R^m.$ Therefore $w \in \LA.$

We now assume that A is associated to a small analytic disc  $f$ 
attached to a generic submanifold $ S \subseteq \C^N.$ 

Take $ U \subseteq \C^m$  and $ V \subseteq \C^{N-m},$
neighborhoods of the origin.  
Denote with $z_1=x +i y,$ the points of $U,$ and with $z_2$ the points of 
$V.$  Assume that the image of $f$ is contained  in $U \times V.$   
Assume moreover the existence of a map  
$ h: \{ |z_1| < \rho , \ |y| < \rho \} \rightarrow \R^m$ with the following 
properties :
$h({\bf 0})=dh_{{\bf 0}}=0,$ and   $x- h(z_2,y)$ is a defining map of 
$S.$ Such map  always exists if we choose $U$ and $V$ small enough.
\cite{Tu}, and \cite{Pa2}.
In such conditions there exists a map $G: \de \Delta \rightarrow 
GL(m,\R)$ such that  $ \gamma \in \EA $ if and only if  $\gamma^t=X^t G$ 
with  $X$ constant vector in  $\R^m$ belonging to the space  
 \[ V_{f}= \{ X \in \R^m : X^t G(h_{z_2} \circ f) \ 
\mbox{extends holomorphically to  $\overline{\Delta}$} \}. \]       
It is proved in  \cite{Pa2}  that the space $V_f$ has  dimension  $d.$
Let Q be a matrix having as rows the vectors of a base of $V_{f}.$  
Then we  can choose 
$P= Q G ( I+i h_y) \circ f $  
If we take a  small enough  $U,$ we can assume that $ ( I+i h_y) \circ f (0)$ 
is invertible, hence $P(0)$ has rank $d.$ \end{proof}

\begin{proposition}

We have $r=l$ if and only if they both equal $d.$
\label{grasm1} \end{proposition}

\begin{proof}
Let $D_x \subseteq \C^d $ be the complex $r$ dimensional image of $P(x).$    
Let
$R_x= D_x \cap (i \R^{d}).$ The set of vectors $c \in 
\C^N$ such that $P(x)\bar{c} \in R_x$ is the orthogonal space in 
$\R^{2N}$ to the rows of the matrix $P(x),$ so it has real codimension $l.$
The Kernel of  $P(x)$ has real codimension $2r,$ hence $R_x$  
has real dimension $2r-l.$ So if $r=l,$ then $R_x$ has real 
dimension $r.$
Now $R_x$ is a totally real subspace of the  $r$ dimensional complex space 
$D_x.$  
It follows that   $D_x =R_x+iR_x.$
Let $G(\C,r,d)$ be the Grasmanian of the complex $r$ subspaces of 
$\C^{d},$ and let  $G(\R,r,d)$ be the analogous real Grasmanian.
The map from $G(\R,r,d)$ sending the subspace  $\Sigma$ to
$\Sigma + i \Sigma,$ identifies $G(\R,r,d)$ 
with a maximal totally real submanifold of  $G(\C,r,d).$ 
Let $\theta: \Delta \rightarrow G(\C,r,d),$ be the holomorphic map 
sending $x$ to $D_x.$  
Our hypotheses implies that the image of $\theta$ is contained in   
$G(\R,r,d).$  Hence $\theta $ is constant. 

In other words, there exist  $r$ independent  vectors  $v_1, \ldots, 
v_{r}$ in 
$\R^{d},$ generating over $\C$ the image of $P(x)$ for every  $x$ in 
$\Delta,$ hence, for every $x$ in $\de \Delta.$  
But the image of $P(x)$ coincides with the image of  $M(x)$ as soon as  $x$ 
is  in $\de \Delta.$  
Therefore, there exist $r$ maps,  $ \lambda_1, \ldots, \lambda_{r}$ 
from $\de \Delta$ to $\R^m,$ such that every map 
$\gamma_j$ of a base of $\EA$ is a linear combination with constant real 
coefficients of the maps $ \lambda_1, \ldots \lambda_{r}.$ Hence 
$d \leq r.$ We conclude by invoking part  a) of Lemma \ref{properties}.
\end{proof}    

\begin{remark}

Since $A$ has rank $m$ in $\de \Delta,$ 
Proposition \ref{ranks} implies that  
$ \theta: \bar{\Delta} \rightarrow G(\C,r,d)$ is an analytic disc attached 
to the maximal totally real submanifold  
$ G(\R,r,d). $ From the above proof we see that, whenever   
$r < d,$ the map  $\theta$ is non-constant. \label{grasm2} \end{remark}

Recall that a  map $A$ is defined in \cite{Trep} to be regular if $N_A=0$ 
i.e. if ${{\Psi}_0}_{| \EA}$ is one to one. 
Then it is natural to give the following 

\begin{definition}
We say that A is strongly regular if the map  ${\Psi}_0$ is one to one. 
(i.e. if $K_A=0 $).
\end{definition}

\begin{examples}

If $ d \leq 1,$ then A is strongly regular.

If m=1, then A is strongly regular if and only if it is regular.
It follows from the results of the next section and of \cite{Trep}.

If A is the map associated to a small analytic disc, then A is strongly 
regular, see part c) of Proposition \ref{properties}.

If m=N, (totally real case) then A is regular if and only if every partial 
index is greater  then or equal to $-1.$ 
\vskip .1cm
A is strongly regular if and 
only if every partial index is greater then or equal to $0.$ 
See Section \ref{tot*real}.
See also \cite{Oh} and \cite{Ce}. \label{strongreg}  \end{examples}

It follows from  Proposition  \ref{grasm1} that $A$ is strongly 
regular if and only if $l=r.$

\section{The Operators} \label{Operator}

Let $U_A : \ONalpha \rightarrow \Calpha,$ be the operator given by  
$U_A(w)= Re(A(w)).$
As we observed in the Introduction, the above operator is 
important for the study of  the set of analytic discs attached to a generic 
submanifold of $\C^N.$   
In particular Tr\'{e}preau gives in $\cite{Trep}$ important properties of the 
extension of $U_A$ to $L^2(\de \Delta,\R^N).$ 

More precisely  
let $H$ be the closure in  $ L^2(\de \Delta, \C) $ of functions in $\ONalpha.$
Then $H$ is the space of functions in $ L^2(\de \Delta,\C) $ whose   
Fouerier series expansion  has only non-negative coefficients.
Let $\tilde{U_A} : H^N \rightarrow L^2(\de \Delta,\R^m)$ be the $L^2$ 
continuous  extension of $U_A.$ 
We will state Tr\'{e}preau result (slightly differently from the statement of  
his paper) in the following 

\begin{theorem}

\[ \begin{array}{l} a) \ \mbox{The operator $U_A : \ONalpha \rightarrow 
\Calpha$ has closed finite codimensional range.} \\ 
\mbox{Moreover its kernel has  a closed supplementary in $\Calpha.$} \\
\mbox{In case $m=N,$  the kernel is  finite dimensional.} \\
\mbox{i.e. $U_A$ is  Fredholm.} \\  
b) \ \mbox{The operator  $\tilde{U_A} : H^N \rightarrow L^2(\de \Delta,\R^m)$
has also closed finite codimensional range.} \\ 
\mbox{Moreover the range $R(U_A)$ coincides with the  intersection  of 
$R(\tilde{U_A})$ with $\Calpha.$} \\
\mbox{In particular the kernel of $\tilde{U_A}$ is in $\Calpha,$ 
hence it coincides with the  kernel of $U_A.$} \\ 
c) \ \mbox{The $L^2$ orthogonal to $R(\tilde{U_A})$ is contained in 
$\Calpha,$} \ \mbox{and it coincides with $N_A.$} \\ 
\ \mbox{It follows that the operator $U_A$ is onto if and only if $A$} \\
\mbox{is regular (i.e. $N_A=0$). } \end{array} \] 
\label{treptrep}  \end{theorem}

\begin{proof}

Tr\'{e}preau deals with the operator $B_A : C^{\alpha}(\de \Delta, \R^N) 
\rightarrow C^{\alpha}(\de \Delta,\R^m)$ given by 
$\phi \rightarrow  Re(A(\phi +i T_0(\phi)))$ where $T_0$ is the Hilbert 
transform normalized at 0. Let 
$ \tilde{B_A}: L^2(\de \Delta,R^n) \rightarrow L^2(\de \Delta,R^m).$ The 
continuous extension of $B_A.$    
Let us observe that  
\[U_A(w)= U_{\frac{A(\zeta)}{\zeta}}(\zeta w)=B_{\zeta^{-1}A}(Re(\zeta w)) \]
\[ Re(\zeta w)(0)= \int_0^{2 \pi} Re(\zeta w) d \theta=0. \]
Since the set 
\[ L= \{ \phi \in L^2(\de \Delta,\R^m) : \int_0^{2 \pi} \phi d \theta=0 \}. \]
is closed m-codimensional, the Theorem follows from the results in 
\cite{Trep}. \end{proof}

Let $U_A(1)$ be the restriction of $U_A$ to the space 
\[ \ONalpha(1)= \{ w \in \ONalpha: w(1)=0 \}. \] 
seen as an operator from $\ONalpha(1)$ to 
\[ \Calpha(1) = \{ \phi \in \Calpha : \phi(1)=0 \}. \]
Since $U_A$ has a closed range with finite codimension, and 
$\ONalpha(1)$ is closed in $\ONalpha$ with finite codimension, if 
follows that $ U_A(1) $ has  closed range with finite codimension as well. 

We introduce 
another operator $V_A $ given by 
\[ V_A(w) = (1-\zeta)^{-1} U_A( (1-\zeta)w ) \]
\[ = 1/2 \left( A(\zeta)w(\zeta)- \overline{\zeta A(\zeta)w(\zeta)} \right). \]
So 
\[ V_A : \ONalpha \rightarrow C^{\alpha}(\de \Delta, \C^m) \]

\begin{lemma}
Let $ w \in \ONalpha$ such that $w(1)=0,$ then the map 
$\frac{w(\zeta)}{1-\zeta}$ is in $ L^2(\de \Delta, \C^N) $ and it is a 
limit in $ L^2(\de \Delta, \C^N), $ of a sequence of maps in $\ONalpha.$
\label{limL2} \end{lemma}

\begin{proof}

Since $\beta > 1/2,$ then 
$\frac{w(\zeta)}{1-\zeta}$ is in $ L^2(\de \Delta, \C^N). $ 
By Lebesque  dominated convergence theorem,  
it is a limit in $L^2,$ of the sequence  
$ w_n= \frac{w}{1-\zeta + 1/n}.$ \end{proof}  

We  take the continuous $L^2$ extension 
$\tilde{V_A}$ of  $V_A.$

\[ \tilde{V_A} : H^N \rightarrow L^2(\de \Delta,\C^m).  \]

\begin{proposition}

Let $\phi \in \Calpha$ such that $\phi(1)=0, $ then $\phi$ belongs to the 
range of $U_A(1)$ if and only if $\frac{\phi(\zeta)}{1- \zeta}$ belongs to 
the range of $\tilde{V_A}.$ \label{range}
\end{proposition}

\begin{proof}
 
Note that, since $\beta > 1/2,$ then for every $\phi$ in $\Calpha,$ with 
$\phi(1)=0,$ we have that   $\frac{\phi(\zeta)}{1- \zeta}$ is in 
$ L^2(\de \Delta,\C^m). $  
If $\phi$ is in the range of $U_A(1),$ then $\frac{\phi(\zeta)}{1-\zeta}$ 
is in the range of $\tilde{V_A},$ because of Lemma \ref{limL2}. 
Viceversa, if  $\frac{\phi(\zeta)}{1-\zeta}$
is in the range of  $\tilde{V_A},$ then there exists 
$u \in  H^N$ such that $\tilde{U_A}((1- \zeta) u ) = 
\phi.$ Set $w= (1- \zeta) u.$
Because of Theorem \ref{treptrep} $ w \in \ONalpha.$ Since 
$ \frac{w(\zeta)}{1- \zeta} = u(\zeta)$ is in   $ L^2(\de \Delta,\C^m),$
we must have $w(1)=0.$ \end{proof}

Note that the range of $\tilde{V_A}$ is contained in the subspace Z of
$ L^2(\de \Delta,\C^m) $ given by the maps $\phi $ such that 
$(1-\zeta) \phi(\zeta)$ is real. So 
\[ Z= \{ \phi \in L^2(\de \Delta,\C^m) : \bar{\phi}= -\zeta \phi \} \] 
Hence  $Z$ is a closed subspace of $ L^2(\de \Delta,\C^m) $ and we can regard 
$\tilde{V_A}$ as an operator from  $ H^N $ into $Z.$

We are going to show that the range of $\tilde{V_A}$ is closed and finite 
codimensional in $Z.$

We first look at the orthogonal to the range of $\tilde{V_A}$ in $Z.$

\begin{proposition} 

An element $\phi \in Z$ is $L^2$ orthogonal to the range of $\tilde{V_A}$ 
if and only if $\bar{\phi}$ is in $Q_A.$ In particular the orthogonal to 
$R(\tilde{V_A})$ is finite dimensional. \label{ort2} \end{proposition}

\begin{proof}

An element $\phi$ in Z is orthogonal to $R(\tilde{V_A})$ if and only if

\[\int_{0}^{2 \pi} \frac{\bar{\phi}^t}{1- \zeta} Re(A(1- \zeta)w) d \theta =0\]
for every $ w \in \ONalpha.$ 
However $\phi= \frac{\gamma}{1- \zeta}$ with real $\gamma,$ so  

\[\int_{0}^{2 \pi} \frac{\bar{\phi}^t}{1- \zeta} Re(A(1- \zeta)w) d \theta \]
\[= \int_{0}^{2 \pi}Re \left(\frac{\gamma^t}{|1- \zeta|^2}(A(1- \zeta)w)\right)d
\theta \]
\[=  Re \int_{0}^{2 \pi} \left( \frac{\gamma^t}{1- \bar{\zeta}}Aw \right) 
d \theta.\]
So we have that $\phi$ in Z is orthogonal to $R(\tilde{V_A})$ if and only if

\begin{equation} \int_{0}^{2 \pi} \bar{\phi^t}A(w)= 0 
 \ \mbox{for every $ w \in \ONalpha$} \label{bar}  \end{equation} 
(replace $w$ with $i w$).
Since $\bar{\phi}= - \zeta \phi,$ then formula (\ref{bar}) is equivalent to
\[ \int_{0}^{2 \pi} \bar{\phi^t}( \zeta A)(w)= 0 \]
and
\[ \int_{0}^{2 \pi} \phi^t (\zeta A)(w)= 0 \]
for every $ w \in \ONalpha.$
Theorem \ref{treptrep} (applied to $ \zeta A $) shows that $\phi \in Z$ is 
orthogonal to $R(\tilde{V_A})$ if and only if   
$\phi$ is in ${\cal C}^{\alpha}(\de \Delta,\C^m),$ 
and $\bar{\phi}^t A$ extends to a map which is 
holomorphic in $\bar{\Delta}$ and vanishes at $0.$ We set 
$\delta=\bar{\phi}.$ By definition
$ \bar{\delta}= -\bar{\zeta}\delta=- \frac{1}{\zeta}\delta.$ Hence 
$ \bar{\delta}^t A$ extends holomorphically to $\bar{\Delta}$ as well,
i.e. $\delta \in Q_A.$
\end{proof}

\begin{lemma}

Given an integer $k,$ the map 
$ p_k : H \rightarrow  L^2(\de \Delta,\C) $
given by 
$ w \rightarrow w + \zeta^k \bar{w},$
has a closed range in $ L^2(\de \Delta,\C).$  
\label{closed} \end{lemma}

\begin{proof}

Let us  identify $ L^2(\de \Delta,\C) $ with $l^2(\C)$ via the orthonormal 
complete system $e^{i m \theta}, \ \mbox{with} $ 
$  m \in \Z.$ 
The map $p_k$ becomes $ a_m \rightarrow a_m + \overline{a_{k-m}}$ under the 
restriction $a_m =0 $ for negative m.
The range of $p_k$ is then  described by the equations $b_{k-m}= \overline{b_m}.$
\end{proof}

We then come to the

\begin{proposition}
The operator $\tilde{V_A}$ has a closed range with finite  codimension in 
$Z.$ \label{Vrange} \end{proposition}

\begin{proof}

Since $ \beta > 1/2 $ we can assume that $A=(A',A")$ where $A'$ has 
values in $GL(m,\C),$ (see \cite{Trep}). So it is sufficient to prove 
the Proposition for the case $m=N.$ In this case we have the factorization 
$-{A}^{-1} \bar{A}= \Theta^{-1} \Lambda \bar{\Theta}.$   (See the proof of 
Lemma \ref{division}).
We can write 
\[ V_A = \frac{1}{2} A \Theta \left( (\Theta^{-1}w)+ (\zeta)^{-1} \Lambda
\overline{(\Theta^{-1}w)} \right). \]
Now multiplication by $A \Theta$ defined an automorphism of 
$ L^2(\de \Delta,\C^m),$  
moreover multiplication by  $\Theta^{-1}$ defines an automorphism of $H^N.$   
Therefore the range of $\tilde{V_A}$ is closed because of Lemma 
\ref{closed} and it is of finite codimension in $Z,$ because of Lemma 
\ref{ort2}. \end{proof}

We summarize the above statements in the following

\begin{theorem}

Let $\phi \in \Calpha$ with $\phi(1)=0.$ Then $\phi$ is in the range of 
$U_A(1)$ if and only if 

\[ \int_{0}^{2 \pi} \gamma^t(\zeta) \frac{\phi(\zeta)}{|1- \zeta|^2} d \theta=0 
\] for  every $\gamma$ in $\EA$ with $\gamma(1)=0.$   
Alternatively if and only if 
\[ \int_{0}^{2 \pi} \delta^t(\zeta) \frac{\phi(\zeta)}{1- \zeta} d \theta=0 
\] for  every $\delta$ in $Q_A.$ \label{Trange} \end{theorem}  

\begin{proof}
It follows from Lemma \ref{range}, Proposition \ref{ort2}, Proposition
\ref{Vrange} and Lemma \ref{division}. \end{proof}

\begin{remark}
Note that the space $\{ (1- \zeta) \bar{\delta} \ \mbox{for $\delta \in Q_A$}
\}$ is a supplementary space  of $R(U_A(1))$ in
$C^{\alpha}(\de \Delta,\C^N)(1).$
Moreover it can be proved in the same way as in \cite{Trep}
that the kernel of $U_A(1)$ has a closed supplementary in 
$C^{\alpha}(\de \Delta,\C^N)(1).$ \label{suppl} \end{remark} 

\begin{corollary}
The operator $U_A(1)$ is onto if and only if $r=d,$ i.e. if and only if
${\Psi}_0$ is one-to-one. \label{onto}  \end{corollary}
\begin{proof}
It follows from Remark \ref{suppl}, and  Lemma \ref{totreal}. \end{proof}

Suppose that $A_{f_0}$ is the map associated to an analytic disc $f_0$ 
attached to a generic submanifold $S$ in $\C^N.$ 
(Recall that $ f \rightarrow \rho \circ f$ is $C^1$). 
Let $q=f_0(1).$ Let 
$M = \{ f \in \ONalpha \ \mbox{such that }$ 
$  \ \mbox{$f$ extends to an analytic disc attached to $ S$} \},$ 
and
$M_q= \{ f \in M \ \mbox{such that $f(1)=q$} \}.$ 

It is proved in \cite{Trep}, by using Theorem \ref{treptrep} together with 
the implicit function theorem, that if $A_{f_0}$ is regular then $M$ is a 
manifold in a neighborhood of $f_0.$ With a similar proof by using corollary 
\ref{onto} and Remark \ref{suppl}, we obtain the following 

\begin{corollary}
If $A_{f_0}$ is strongly regular, then $M_q$ is a manifold in a neighborhood of 
$f_0.$ \label{varieta'} \end{corollary}

\begin{remark}
In the totally real case $M_q$ is finite dimensional and its dimension
equal the total index k of the submanifold S. 
In general (when S is totally real but A may not be strongly regular) 
the number k is the index of $U_A(1).$ (See section \ref{tot*real}). 
\end{remark}

\begin{remark}
More generally $M$ is a  manifold in a neighborhood of $f_0$ as 
long as $d-l$ is  constant as $f$ varies in a 
neighborhood of $f_0.$ 
Similarly $M_q$ is a  manifold in a neighborhood of $f_0$ as 
long as $d-r$ is  constant as $f$ varies in a neighborhood of $f_0.$ 
\label{implfunc} \end{remark}

We finally  need to study the operator $U_A$ restricted  to the subspace 
of $\ONalpha$ given by the maps vanishing at the point $1$ and at the point 
$-1$ simultaneously.

Denote by $\ONalpha(1,-1)$ the space of maps in $\ONalpha$ which vanish at 
$1$ and at $-1.$ Similarly denote by $C^{\alpha}(1,-1),$ the space of 
$C^{\alpha}$ maps vanishing at $1$ and $-1.$
Denote by $U_A(1,-1)$ the restriction of  $U_A$ to the space   
$\ONalpha(1,-1)$
seen as an operator from $\ONalpha(1,-1)$ to $C^{\alpha}(1,-1).$

Note that if a map  $ \phi$ is in $C^{\alpha}$ and it vanishes at $1$ and 
$-1,$ then $\frac{\phi(\zeta)}{1- \zeta^2} \in L^2.$ 
Define

\[ {V^*}_A(w) = (1-\zeta^2)^{-1} U_A( (1-\zeta^2)w ) \]
\[= 1/2 \left( A(\zeta)w(\zeta)- \overline{\zeta^2 A(\zeta)w(\zeta)} \right).\]
So 
\[ {V^*}_A : \ONalpha \rightarrow C^{\alpha}(\de \Delta \C^m) \]
and 
\[ \tilde{{V^*}_A} : H^N \rightarrow L^2(\de \Delta,\C^m).  \]

\begin{lemma}
Let $ w \in \ONalpha$ such that $w(1)=0$ and $w(-1)=0,$ then the map 
$\frac{w(\zeta)}{1-\zeta^2}$ is in $ L^2(\de \Delta, \C^N) $ and it is a 
limit in $ L^2(\de \Delta, \C^N), $ of a sequence of maps in $\ONalpha.$
\label{limL2bis} \end{lemma}

\begin{proof}
The proof is  analogous to that of Lemma \ref{limL2}. \end{proof}

\begin{proposition}

Let $\phi \in \Calpha$ such that $\phi(1)=0, $ and $\phi(-1)=0,$ 
then $\phi$ belongs to the 
range of $U_A(1,-1)$ if and only if $\frac{\phi(\zeta)}{1- \zeta^2}$ 
belongs to the range of $\tilde{{V^*}_A}.$ \label{range,bis}
\end{proposition}

The range of $\tilde{{V^*}_A}$ is contained in the subspace $Z^*$ of
$ L^2(\de \Delta,\C^m) $ given by the maps $\phi $ such that 
$(1-\zeta^2) \phi(\zeta)$ is real. So 
\[ Z^*= \{ \phi \in L^2(\de \Delta,\C^m) : \bar{\phi}= -\zeta^2 \phi \} \] 

$Z^*$ is a closed subspace of $ L^2(\de \Delta,\C^m) $ and we can regard 
$\tilde{{V^*}_A}$ as an operator from  $ H^N $ into $Z^*.$

Let $K^*_A$ be the subset of $\EAC$ of maps $\delta$ such that the 
holomorphic extension of $\delta^t A $ vanishes of order at least 2 at $0.$
Let   
$ Q^*_A \subseteq K^*_A$ the subspace given by
\[ Q^*_A = \{ \delta \in K^*_A : \bar{\delta} = 
-\bar{\zeta^2} \delta. \} \]

\begin{proposition} 

An element $\phi \in Z^*$ is $L^2$ orthogonal to the range of $\tilde{V^*_A}$ 
if and only if $\bar{\phi}$ is in $Q^*_{(\zeta A)}.$ 
In particular the orthogonal to $R(\tilde{V^*_A})$ is finite dimensional. 
\label{ort2,bis} \end{proposition}

\begin{proposition}
The operator $\tilde{V^*_A}$ has a closed range with finite codimension in 
$Z^*.$ \label{Vrange,bis} \end{proposition}

\begin{lemma}

The space  $Q^*_{(\zeta A)}$ is maximal totally real in $K^*_{(\zeta A)}.$ 
\label{totreal,bis} \end{lemma}

However, the space  $Q^*_{(\zeta A)}$ is essentially $\EA,$ more precisely 
we have

\begin{lemma}
The map from $\EA$ into  $Q^*_{(\zeta A)}$ given by
\[ \gamma \rightarrow i \zeta \gamma \] is an isomorphism,
which extends to an isomorphism of $\EAC$ onto $K^*_{(\zeta A)}.$
\label{defect} \end{lemma}

\begin{proof}

Just observe that $\delta= i \zeta \gamma$ with real $\gamma,$ if and only 
if $ \bar{\delta}=-\bar{\zeta^2} \delta.$ Moreover $ i \zeta \gamma^t 
(\zeta A)$ extends to a holomorphic map vanishing of order at least 2 at 
$0,$ if and only if $\gamma^t A$ extends holomorphically. 
The last assertion of the Lemma follows from Lemma \ref{totreal,bis}.
\end{proof}

From the above Lemma, with a similar proof as in Theorem \ref{Trange} we 
find 

\begin{theorem}

Let $\phi \in \Calpha$ with $\phi(1)=\phi(-1)=0.$ Then $\phi$ is in the range 
of $U_A(1,-1)$ if and only if 
\[ \int_{0}^{2 \pi} \gamma^t(\zeta)\frac{ \zeta \phi(\zeta)}{1- \zeta^2} 
d \theta=0 \] for  every $\gamma$ in $\EA.$ \label{Trange,bis} \end{theorem}  

\begin{remark}
Note that the space $\{ i \bar{\zeta}(1- \zeta^2) \gamma \ \mbox{for $\gamma 
\in \EA$} \}$ is a supplementary space  of $R(U_A(1,-1))$ in $C^\alpha(1,-1).$
\label{suppl,bis} \end{remark} 

\section{Evaluation spaces } \label{eval}

Fix an element $y \in \de \Delta,$ and an element $ x \in \bar{\Delta},$
define 
\[ \LA(y)= \{ w \in \LA : w(y)=0 \} .\]
We are interested in studying the various evaluation spaces
$ W_x= \mu_x(\LA),$ and $ W_x(y)= \mu_x(\LA(y)).$ 
(Recall that $\mu_x$ is the evaluation at $x$ of an element in 
$\ONalpha$ and $\LA= \{ w  \in \ONalpha  : Re(Aw)= 0 \}.$ )
We have the following theorems:

\begin{theorem}
Given $ x \in \Delta$ we have 
\[a) \  W_x = \{ c \in \C^N : Re(\Psi_x(\gamma)c)=0 \ \mbox{for every 
$\gamma \in \EA$} \}. \] 
In other words $W_x= \{ c \in \C^N : Re(P(x)c)=0 \}.$
So   $W_x \cap i W_x$ coincides  with the kernel of $P(x).$

b) The real codimension of $W_x$ in $\R^{2N}$  is $l.$

c) The complex codimension of   $W_x \cap i W_x$ in $\C^N$ is $r.$

d) The complex codimension of   $W_x + i W_x$ in $\C^N$ is $ l-r.$
\label{val0}  \end{theorem}

\begin{proof}
We assume first that $x=0.$
A vector   c in $\C^N,$  
belong to $W_0$ if and only if there exists 
$u \in {{\cal O}^N}_{\alpha}$ 
such that $\zeta u(\zeta) + c $ belongs to $\LA,$ if and only if  
$-Re(A(\zeta) c)$ belongs to the range of  $U_{ \zeta A}.$ 
Now $ \gamma \in N_{ \zeta A}$ if and only if $ \zeta \gamma^t A$ extends 
holomorphically to $\bar{\Delta}$ and its extension vanishes at 0. This is 
the case if and only if $\gamma \in \EA.$ It follows from Theorem 
\ref{treptrep} that  $c \in  W_0$ if and only if 
\[ \int_0^{2 \pi} \gamma^t Re(A c) d \theta = 0 \] 
for every $\gamma \in \EA.$
However, since $ \gamma \in \EA,$ we have  
\[ \int_0^{2 \pi} \gamma^t Re(A c)d \theta =Re(\Psi_0(\gamma)c). \]
In other words the image of $W_0$ under the conjugation map, is the 
orthogonal space in $\R^{2N}$ to the space $\Psi_0(\EA).$ This proves b).
Part c) and part d) follow easily. 
In the general case, let $\sigma$ be an authomorphism of the disc such that 
$\sigma(0)=x.$ It follows from Proposition \ref{ranks} that 
the numbers $r$ and $l$ relative to the map $A,$ coincide with the 
corresponding numbers relative to $A \circ \sigma.$ \end{proof}

\begin{theorem}
Given $ x \in \de \Delta$ we have 
\[ a) \  W_x = \{ c \in \C^N : Re(A(x)c)=0 \ \mbox{and} \ 
\delta^t(x)A(x)c=0 \ \mbox{ for every $\delta \in K_A $ } \}. \] 
 
So the kernel of $A(x)$ is contained in $ W_x.$

b) The real codimension of $W_x$  in the space $Re(A(x)c)=0$ is 
$l-r.$ 

c) We have $W_x \cap i W_x = Ker(A(x)).$

d) The complex codimension of   $ W_x + i W_x$ in $\C^N$ is $ l-r.$
\label{val1} \end{theorem}

\begin{proof}

Assume first that $x=1.$

Obviously  
\[ W_1  \subseteq  \{ c \in \C^N : Re(A(1)c)=0 \}. \]
Now, given $c \in \C^N \ \mbox{with} \  Re(A(1)c)=0,$ we have $ c \in W_1 $ 
if and only if $- Re(A(\zeta)c)$ belongs to the range of $U_A(1).$
This is the case if and only if 
\[ \int_0^{2 \pi} \frac{\delta^t Re(A(\zeta)c)}{1- \zeta} d \theta =0 \ 
\mbox{for every $ \delta \in Q_A.$} \]  
By the definition of $Q_A$ this is equivalent to saying that
\[ \int_0^{2 \pi} Re \left( \frac{\delta^t A(\zeta)c}{1- \zeta}\right)d \theta 
=0 \ \mbox{for every $ \  \delta$  in $ Q_A.$} \]
Let $f$ be the holomorphic extension of $\delta^t A c.$
Since $\delta= - \zeta \bar{\delta},$ then $\delta(1)$ is purely 
imaginary, and, by assumption, so is $A(1)c,$ therefore $f(1)$ is real.
Moreover, since $\delta \in Q_A \subseteq K_A,$ then $f(0)=0.$
We have: 
\[ \frac{1}{2 \pi} \int_0^{2 \pi} Re \left( \frac{f(\zeta)}{1- \zeta} 
\right)  d \theta = \]
\[\frac{1}{2 \pi} \int_0^{2 \pi} Re \left( \frac{f(\zeta)- f(1)}{1- \zeta}  
\right) d \theta \]
\[+ \frac{f(1)}{2 \pi} \int_0^{2 \pi} Re \left( \frac{1}{1- \zeta} \right) 
d \theta. \]
By the proof of Lemma \ref{limL2} we know that    
\[\frac{1}{2 \pi} \int_0^{2 \pi} Re \left( \frac{f(\zeta)- f(1)}{1- \zeta}  
\right) d \theta \]
\[ = \lim_{n \rightarrow \infty} \frac{1}{2 \pi} \int_0^{2 \pi} 
Re \left( \frac{f(\zeta)- f(1)}{1- \zeta +1/n}  \right) d \theta \]
\[= \lim_{n  \rightarrow \infty} \frac{Re(f(0)- f(1))}{1 +1/n}= -f(1), \]  
since $f(0)=0$ and $f(1)$ is real.
On the other hand
\[ \frac{f(1)}{2 \pi} \int_0^{2 \pi} Re \left( \frac{1}{1- \zeta} \right) 
d \theta = \frac{f(1)}{2} \] 
We have then proved that 
\[ W_1 = \{ c \in \C^N : Re(A(1)c)=0 \ \mbox{and} \  \delta^t(1)A(1)c=0 \] 
\[ \mbox{ for every $ \delta \in Q_A $} \}. \] 
To replace $Q_A$ with $K_A$ we only need to recall that $Q_A$ is maximal totally
real in $K_A.$
As we already observed if $\delta \in Q_A,$ then $i \delta(1)$ is real. 
Let $X \subseteq \R^m,$ 
be the orthogonal space to the image of the map 
$\delta \rightarrow i \delta(1)$ for $\delta \in Q_A.$ Then because of 
Remark \ref{division2},  
$X$ has real codimension $d-r - (\mbox{dimension of $N_A$}) = l-r.$ 
Moreover $W_1={A(1)}^{-1}(iX).$ Since the kernel of 
$A(1)$ has complex dimension $N-m,$ then $W_1$ has real dimension 
$2N-m-l+r,$ therefore it has codimension $l- r $ in 
$Ker(Re A(1)).$ This proves part b). Part c) and d) are immediate. 
In the general case, let $\sigma$ be the rotation such that $\sigma(1)=x.$
Since $\sigma(0)=0,$ then $ \delta \in K_A$ if and only if 
$\delta \circ \sigma \in K_{A \circ \sigma}.$ The conclusion follows.  
\end{proof}

\begin{remark}
Let  $ x \in \de \Delta,$ and  $\gamma$ in $\EA,$ such that $\gamma(x)=0.$ 
Define    
\[\gamma'(x)= \lim_{y \rightarrow x} \frac{\gamma(y)}{y-x}. \] 
(Such limit exists and its finite because of Lemma \ref{division}).
Then part a) of Theorem \ref{val1} can be restated as follows:
\[ W_x = \{ c \in \C^N : Re(A(x)c)=0 \ \mbox{and} \]   
\[ \overline{\gamma'(x)}A(x)c=0 \mbox{ for every $\gamma  \in \EA $ 
such that $\gamma(x)=0 $} \}. \] 
This follows from the proof of Theorem \ref{val1}, and from Lemma 
\ref{division}. \end{remark}

For $x \in \Delta,$  let 
\[ \begin{array}{l}  \tilde{\Psi}_x : 
{\cal E}_{(\zeta-x)A} \rightarrow \C^N \\
\mbox{be the evaluation at $x$ of the holomorphic extension of 
$\delta^t(\zeta-x)A,$} \end{array} \]
In other words $ \tilde{\Psi}_x(\delta) $ is the residue at $x$ of the 
meromorphic extension of $\delta^t A.$

Moreover let $r_1$ and $l_1$ and $d_1$ be the numbers $r$ and $l$ and $d$ 
relative to the map $\zeta A.$

\begin{theorem}

Given $ x \in \Delta$ and $y \in \de  \Delta,$  we have 

\[ a) \ W_x(y) = \{ c \in \C^N : Re(\tilde{\Psi}_x(\gamma)c)=0 \]
\[ \mbox{for every} \  \gamma \in {\cal E}_{(\zeta-x)A}  \ \mbox{such that} \  
\gamma(y)=0 \}. \] 
 
\[ W_x(y) \cap  i W_x(y) = \{ c \in \C^N : \tilde{\Psi}_x(\gamma)c=
0 \] 
\[ \mbox{for every} \  \gamma \in {\cal E^{\C}}_{(\zeta-x)A} \ \mbox{such that} 
 \ \gamma(y)=0  \}. \] 

b) The real codimension of $W_x(y)$ in $\R^{2N}$ is 
$l_1 - r_1 + r.$

c) The complex codimension of   $W_x(y) \cap i W_x(y)$ in $\C^N$ is  
$l_1 - r_1.$

d) We have $W_x(y) + i W_x(y)=W_x \cap i W_x = Ker(P(x)).$ 
\label{val2}  \end{theorem}

\begin{proof}
Assume first that $x=0$ and $y=1.$
Because of Lemma \ref{division}, statement a) is equivalent to 
\[ W_0(1) = \{ c \in \C^N : Re(\tilde{\Psi}_0(\bar{\delta})c))=0 \] 
for every $ \delta \in Q_{\zeta A} \}. $
\vskip .1cm
A vector $c$ belongs to  $W_0(1)$  if and only if there exists a map 
$ w \in \LA$ with $w(1)= 0 $ and $w(0)=c.$ This is the case if and only 
if   there exists $u$  in $\ONalpha$   such that, 
\[ Re( A(\zeta u(\zeta)+c))= Re((\zeta A)(u(\zeta)+c))+ 
Re(A((1-\zeta)c))=0. \] 
Since $w= \zeta u +c,$ we have $u(1)+c=0.$ 
Therefore $c$ is in $W_0(1)$ if and only if $-Re(A((1 -\zeta)c))$ 
belongs to the range of $U_{\zeta A}(1).$ 
This is so if and only if 
\[ \int_0^{2 \pi} \delta^t \frac{Re(A((1 -\zeta)c))}{1- \zeta} d \theta= 0 \]
for every $\delta \in Q_{(\zeta A)}.$ 
Since  every $\delta \in Q_{(\zeta A)}$ is of the form 
$\delta= \frac{\gamma}{1- \bar{\zeta}}$ with $\gamma$ real, we find 
\[ \int_0^{2 \pi} \delta^t \frac{Re(A((1 -\zeta)c))}{1- \zeta} d \theta \]
\[=  Re (  \int_0^{2 \pi} \delta^t Ac d \theta) =  \]
\[  Re \left(  \int_0^{2 \pi} -\bar{\delta}^t(\zeta A)c d \theta \right)   \]
\[=-Re(\tilde{\Psi}_0(\bar{\delta})c).\] 
Therefore the image of   $W_0(1)$ under the conjugation map, is the 
orthogonal space in $\R^{2N}$ to $H(Q_{(\zeta A)})$ where  
$ H: K_{\zeta A} \rightarrow \C^N $ is given by
\[ H(\delta)= \tilde{\Psi}_0(\bar{\delta}). \] 
It follows directly from the definitions that the kernel of H is precisely 
$\EAC,$ and its intersection with $Q_{\zeta A} $ is 
$Q_A.$ Now $Q_{\zeta A}$ has dimension 
$ d_1- r_1,$ and $Q_A$ has dimension $d-r.$ 
Hence the image of $H$ has dimension 
$ d_1- d- r_1 +r.$ However $d$ is the dimension of 
$N_{\zeta A}$ which is $  d_1- l_1.$ 
So the codimension of
$W_0(1)$ is $l_1- r_1 + r.$ 
\vskip .1cm
Similarly, since $Q_A$ is maximal totally real in $K_A,$ then  
a vector $c \in \C^N,$ belongs to $W_0(1) \cap i W_0(1)$ 
if and only if $\tilde{\Psi}_0(\bar{\delta})c=0$ for every $\delta \in 
K_A.$ Therefore the complex codimension of $W_0(1) \cap i 
W_0(1)$ equals the complex dimension of $K_{\zeta A}$ minus the 
complex dimension of $\EAC,$ which is  $l_1 - r_1.$
From above we derive that  
$W_0(1) + i W_0(1))$
has complex codimension $r.$ It follows from Remark \ref{pearing}  
and Theorem \ref{val0} that 
$W_0(1)+i W_0(1) \subseteq  W_0 \cap i W_0.$ 
\vskip .1cm
Since, by Theorem \ref{val0}, these two 
complex spaces have the same dimension, they must coincide. 

For the general case, let $\sigma$ be the unique authomorphism of the disc 
such that $\sigma(1)=y $ and $\sigma(0)=x.$
By replacing $A$ with $A \circ  \sigma$ we immediately prove part c).
\vskip .1cm

Now $\gamma \circ \sigma \in {\cal E}_{ \zeta A \circ \sigma}$
if and only if  $\gamma \in {\cal E}_{ \sigma^{-1} A}.$
On the other hand  $\sigma^{-1}(x)=0,$ and since $\sigma$ is an authomorphism,
the function $\eta(\zeta)$ given by 
\[ \eta(\zeta) = \frac{\sigma^{-1}(\zeta)}{\zeta - x} \]
for $\zeta \neq x$ and 
\[ \eta(x)= \frac{\de \sigma^{-1}}{\de \zeta}_{| \zeta=x } \] 
is holomorphic on $\Delta,$ $C^{\alpha}$ on $\de \Delta$ and does not have 
zeros on $\bar{\Delta}.$
Therefore  $ {\cal E}_{\sigma^{-1} A}= {\cal E}_{(\zeta-x) A}.$
\vskip .1cm
We are left to show that the dimension of $W_x(y)$ is independent on 
$x$ and $y.$ First of all the space 
\[ \{ \gamma \in {\cal E}_{(\zeta-x)A} : \gamma(y)=0 \} \] 
has dimension 
independent on $x$. In fact the map 
$ \gamma \rightarrow t \gamma,$ (where the function 
\[ t(\zeta)=\frac{\bar{x}\zeta^2 -(|x|^2+1)\zeta +x}{\zeta}, \]
was defined in the proof of Proposition \ref{ranks}), gives a linear 
isomorphism between 

$ \{ \gamma \in {\cal E}_{(\zeta-x)A} : \gamma(y)=0  \} $
and
$ \{ \gamma \in {\cal E}_{\zeta A} : \gamma(y)=0 \}. $

However, since for every rotation $\sigma,$ we have that 
${\cal E}_{\zeta A}$ is isomorphic to $ {\cal E}_{\zeta A \circ \sigma},$
by choosing the rotation sending $1$ to $y$ we show that the dimension of 
the above spaces is independent on $y$ as well.
To conclude we still need to prove the independence on $x$ and $y$ of the 
dimension of the spaces  
\[ \{ \gamma \in {\cal E}_{(\zeta-x)A} : \gamma(y)= \tilde{\Psi}_x(\gamma)
=0  \} \]

This space coincides with 
$ \{ \gamma \in {\cal E}_{A} : \gamma(y)=0 \} $
whose dimension is independent on $y$ because of part b) of Proposition 
\ref{ranks}. \end{proof}

\begin{theorem}
Given $ x \ \mbox{and} \  y \  \in \de \Delta \ \mbox{with $x \neq y$} $ 
we have 
\[ a) \   W_x(y) = \{ c \in \C^N : Re(A(x)c)=0 \  \mbox{and} \   
\gamma^t(x)A(x)c=0 \mbox{ for every $\gamma \in \EA$ } \}. \] 
 
So the kernel of $A(x)$ is contained in $ W_x.$

b) The real codimension of $W_x(y)$  in the space $Re(A(x)c)=0$ is 
$r$

c) We have $W_x(y) \cap i W_x(y) = Ker(A(x)).$

d) We have $ W_x(y) + i W_x(y)= Ker P(x).$ 
\label{val3}  \end{theorem}

\begin{proof}
We first assume that $y=1$ and $x=-1.$ Given $c \in \C^N$ with 
$Re(A(-1)c)=0,$ we have that $c$ belongs to $W_{-1}(1)$ if and only if
there exists $w \in \ONalpha$ such that $w(1)=0, \ w(-1)=c$ and 
$Re(A(\zeta)w(\zeta))=0.$ Let us write $w= u + \frac{(1-\zeta)c}{2}.$
Then  $c$ belongs to $W_{-1}(1)$ if and only if
$Re((1-\zeta)A(\zeta)c)$ belongs to the range of $U_A(1,-1).$ 
This is to say that 
\[ \int_0^{2 \pi} \zeta \gamma^t \frac{Re(A((1 -\zeta)c))}{1- \zeta^2} 
d \theta=0 \ \mbox{for every}  \ \gamma \in \EA.  \]
However
\[ \int_0^{2 \pi} \zeta \gamma^t \frac{Re(A((1 -\zeta)c))}{1- \zeta^2} 
d \theta =\]
\[ \frac{1}{2} \int_0^{2 \pi}  \zeta \gamma^t \frac{Ac}{1+\zeta} d \theta  + \]
\[ \frac{1}{2} \int_0^{2 \pi} \zeta \gamma^t\frac{(1-\bar{\zeta})\bar{A}\bar{c}}
{1-\zeta^2} d \theta = \]
\[ \frac{1}{2}  \int_0^{2 \pi}  \zeta \gamma^t \frac{Ac}{1+\zeta} d \theta \]
\[-\frac{1}{2} \int_0^{2 \pi} \gamma^t \frac{\bar{A} \bar{c}}{1+\zeta} d 
\theta = \]
\[=i Im \left(\int_0^{2 \pi}  \zeta \gamma^t \frac{Ac}{1+\zeta} d \theta 
\right). \] 
Let $f$ be the holomorphic extension of $\zeta \gamma^t A c,$ then $f(0)=0$ 
and $Re(f(-1))=0.$  Then
\[i Im \left( \int_0^{2 \pi} \frac{f(\zeta)}{1+ \zeta} d \theta \right) =\]
\[i Im \left( \int_0^{2 \pi} \frac{f(\zeta) - f(-1)}{1+ \zeta} d \theta 
\right) \]
\[+ f(-1) Re \left( \int_0^{2 \pi} \frac{1}{1+ \zeta} d \theta \right)= 
- \pi f(-1). \] 
Since $A(-1)$ has maximal  rank, we find  from part b) of 
Proposition \ref{ranks}, that the space $ \{ \gamma(-1), \gamma \in \EA \} 
\subseteq \R^m $ has real dimension $r$
and its orthogonal space $X$ in $\R^m,$ has real dimension $m-r.$    
We conclude that 
$ W_{-1}(1)= i A(-1)^{-1}(X)$  has real dimension $2N-m -r,$ so 
it has codimension $r$ in $ Ker (Re(A(-1)),$ this proves part b). 

Part c) follows from a). As far as d) concernes we conclude from a) and c) 
that   
$ W_{-1}(1)+iW_{-1}(1)$ has complex codimension $r$ in $\C^N.$ 
However form Remark \ref{pearing} we know that  
$ W_{-1}(1)+iW_{-1}(1)$ is a subspace of $Ker P(-1)$ which also has codimension 
$r$ in $\C^N.$ 
The general case follows by acting with the authomorphism group of the 
disc. \end{proof}

\begin{proposition}

a) The assignment $x \rightarrow Ker(P(x))$ defines a complex $C^{\alpha}$ 
vector bundle $F$ on $\bar{\Delta}$ of rank $N-r,$ which is holomorphic on 
$\Delta.$ Moreover the bundle $Ker A$ on $\de \Delta$ is a subbundle of the 
restriction of $F$ to  $\de \Delta.$  

b) The assignment $x \rightarrow W_x \cap i W_x$ coincides with 
$F$ on  $\Delta.$
  
c) Given $y \in \de \Delta,$ the assignment 
$x \rightarrow W_x(y) + i W_x(y)$ coincides with 
$F$ on  $\bar{\Delta}-y.$
\end{proposition} 

\begin{proof}

Part a) follows from the definition of $P$ and from Proposition 
\ref{ranks}.

Part b) follows from part a) of Theorem \ref{val0}.

Part c) follows from Theorem \ref{val3}. \end{proof}

There is an interesting special case

\begin{proposition}
If for the map $A$ we have $r=m,$ then $W_x(y)$ is a complex space 
for every $y \in \de \Delta,$ and  every $x \in \bar{\Delta}-y.$ 
Therefore, in this case, the bundle $F$ is a holomorphic extension to 
$\bar{\Delta}$ of the bundle $Ker(A(x))$ on $\de \Delta.$
\label{r=m} \end{proposition}

\begin{proof}
If $x \in \de \Delta,$ the result follows from Theorem \ref{val3} part c).
Assume then, that $x \in \Delta,$ by acting with the authomorphism group of 
the disc we may suppose $x=0,$ $y=1.$
We want to show that $W_0(1)= Ker(P(0)).$
The inclusion $W_0(1) \subseteq  Ker(P(0))$ follows from Remark 
\ref{pearing}. Fix  $c \in Ker P(0).$  
Theorem \ref{val0} part a) says that there exists  $w \in \LA $ with $w(0)=c.$
Hence  $Pw$ is a purely imaginary constant such that  $P(0)w(0)=0,$
so $Pw \equiv 0.$ Set  $\tilde{w}= (1- \zeta) w,$ then $\tilde{w}(1)=0$ 
and $\tilde{w}(0)=c,$ from Proposition \ref{properties} part b) we conclude 
that $\tilde{w}$ belongs to $\LA.$  \end{proof} 

Suppose now that $A_{f_0}$ is the map associated to an analytic disc $f_0$ 
attached to a generic submanifold $S$ in $\C^N.$ Let $q=f_0(1),$ 
Recall that, whenever $A$ is regular, the set 
$M = \{ f \in \ONalpha \ \mbox{such that }$ 
$  \ \mbox{$f$ extends to an analytic disc attached to $ S $} \},$ 
is a manifold in a neighborhood of $f_0,$ and that, whenever $A_{f_0}$ 
is strongly regular, the set
$M_q= \{ f \in M \ \mbox{such that $f(1)=q $} \}, $ 
is a manifold in a neighborhood of $f_0.$ 

We then have the following 

\begin{theorem}

a) If $ \zeta A$ is regular, that is if $\EA=0,$  then there exists a 
neighborhood  $U$  of $f_0$ in $M$ such that  ${\mu}_0(U)$ is an open set in 
$\C^N.$   

b) If $ \zeta A$ is strongly regular, then there exists a 
neighborhood  $U$  of $f_0$ in $M_q$ such that  ${\mu}_0(U)$ is an open set in 
$\C^N.$ \label{openmap}  \end{theorem}

\begin{proof}
 
a) If $\zeta A$ is regular, then  the operator $U_{\zeta A}$ is   
onto, in particular so is the operator   $U_{A}.$ Hence  
$M$ is a manifold in a neighborhood of  $f_0.$ 
The tangent space to $M$ at $f_0$ coincides with $\LA$ and the map 
${{\mu}_0}_{| \LA}$ coincide with the differential of  ${{\mu}_0}_{|M}.$
It follows from Theorem \ref{val0} that such differential is onto.
We conclude by using the open mapping theorem. 
The proof of part b) goes in the same way if we use Theorem \ref{onto} 
\end{proof}

\begin{lemma}

Let $\Omega$ be the subset of $\ONalpha$  given by the analytic discs $f$ such 
that there exists a neighborhood of f in $\ONalpha$ where the numbers  
d, l and r are constant. Then $\Omega$ is open and dense in 
$\ONalpha,$ moreover $\Omega \cap M$ is open and dense in $M.$  
\label{locconst} \end{lemma}

\begin{proof}
 
Let $f$ be an analytic disc with associated map A. Recall that $d-l$ is the 
dimension of the cokernel of $U_A,$ whereas $d-r$ is the dimension of the 
cokernel of $U_A(1)$ and d is the dimension of the cokernel of $U_A(1,-1).$
Since the map $f \rightarrow A_f$ is continuous on $\ONalpha,$ 
It follows from, \cite{K} pag.235,  that $d,$ $d-r,$ and $d-l$ are upper 
semicontinuous functions in $\ONalpha$ with value in the set of non negative 
integers. Let $ \Omega_1 $ be  the set of discs where $d$ is locally constant,
$\Omega_2$  the set of discs where $d-l$ is locally costant and 
$ \Omega_3$  the set of discs where $d-r$ is locally constant. 
Let us consider for example $\Omega_1.$ Let U be any non empty open subset 
of $\ONalpha$, then, since d gives a semicontinuous function, it follows that 
the set of discs in U where d takes its minimum value in U, is open and non 
empty, so $U \cap \Omega_1 \neq \emptyset.$ In the same way it can be shown 
that $\Omega_2$ and $\Omega_3,$ are dense in M. So 
$\Omega= \Omega_1 \cap \Omega_2 \cap \Omega_3$ 
is a dense open subset of $\ONalpha,$ in the same way it can be proved that 
$\Omega \cap M$ is dense in $M.$  \end{proof}

Observe that $\Omega$ is invariant under the action of the group of 
authomorphisms of the disc.

\begin{definition}

Let $p$ be a point on the submanifold S in $\C^N,$ we say that $p$ is a 
minimal point if every immersed  submanifold H  of $S$ containing $p$ and such 
that the complex tangent bundle to $H$ concides with the restriction to H of 
the complex tangent bundle to S,  is an  open set in S. \label{min}  
\end{definition}

\begin{corollary}
Let $S$ be a generic submanifold of $\C^N$ where every point is minimal, 
then the set $\tilde{\Omega}$ of analytic discs attached to S with defect zero 
is open and dense in $M.$  In particular $\mu_0$ restricted to 
$\tilde{\Omega}$ is an open map. \label{opendense} \end{corollary}  

\begin{proof}
Let $\Omega$ be the open dense subset defined in Lemma \ref{locconst}, 
clearly $\tilde{\Omega}, $ is contained in $\Omega \cap M.$ Viceversa let 
$f_0 \in \Omega \cap M$ such that $f(0)=q,$   then by Remark \ref{implfunc}, 
the sets $\Omega \cap M$ 
and $\Omega \cap M_q$ are manifolds. Let us consider  the 
the map $\mu_{-1}$ restricted to $\Omega \cap M_q,$ 
seen as a map from $ \Omega \cap M_q$ into S. By Theorem \ref{val3} the 
image of the differential of such map has  a locally constant 
codimension r. Therefore $ \mu_{-1}(\Omega \cap M_q) $  is an immersed 
submanifold as in definition \ref{min}. Since every point of $S$ is assumed 
to be minimal, it follows that $ r=0,$ hence, because of Proposition 
\ref{properties} d=0. \end{proof}

\section{The case $N=m$ (totally real case)} \label{tot*real}

If we assume that $A$ is at values in $GL(N,\C),$ then we may use the Birkhoff
factorization $-A^{-1} \bar{A}= \Theta \Lambda \overline{\Theta^{-1}},$ 
described in Lemma \ref{division}.
We are going to determine everything we need  in terms of the partial indices 
$k_1, k_2, \ldots k_N.$ (See also \cite{Pa}).

We have 
\[ \LA= \{ u \in \ONalpha : Au + \overline{A} \overline{u}=0 \} \]
\[ = \{ u \in \ONalpha : u + A^{-1} \overline{A} \overline{u}=0 \} \]
\[ = \{ u \in \ONalpha : \Theta^{-1} u - \Lambda \overline{ 
\Theta^{-1}u}=0 \} \]
\[ = \{ v \in \ONalpha : \bar{v_j}= \zeta^{-k_j}v_j \ \mbox{on $\de \Delta $} 
\} \]
Here $j$ runs  from $1$ to $N$ and   $v_j$ is the jth  component of the vector 
$v= \Theta^{-1}u.$
Therefore if $k_j < 0,$ then $v_j \equiv 0,$ if $k_j=0$ then $v_j$ is a 
real constant function, and if $j >0 $ then $v_j$ is a polynomial  of degree 
$k_j.$ It has the  form $v_j= \sum_{i=0}^{k_j} a_i z^i,$ with 
$a_s= \overline{a_{k_j-s}}$ for every $s$ from $0$ to $k_j.$ 

Let $k_{+}$ be the sum of all the positive partial indeces, and let  $m_{+}$ be 
their number. Let  $k_{-}$ be the absolute value of the sum of the negative 
partial indices,  and let  $m_{-}$ their number, finally let $m_0$ be the 
number of partial indices equal to $0.$
So we have $dim_{\R} \LA= k_{+} + m_{+} + m_{0}.$  
On the other hand 
\[ \EA= \{ \gamma \in \Calpha: \gamma^t A=u \ \mbox{extends 
holomorphically}\}.\]
Hence $\EA$ is isomorphic to the space 
\[ \{ u \in \ONalpha : u^t A^{-1} \ \mbox{is real} \} \]
\[= \{ u \in \ONalpha : - u^t \Theta \Lambda \overline{\Theta^{-1}} = 
\bar{u}^t.\} \]
\[= \{ v \in \ONalpha :  \bar{v_j}= \zeta^{k_j}v_j \} \]
where $j$ runs  from $1$ to $N$ and   $v_j$ is the jth  component of the vector 
$v= i (\Theta^{-1})^t u.$
In the same way as above we find that  $d= k_{-} \  + \  m_{-} \  + \  m_0.$
\vskip .1cm
If we replace A with $\zeta^{-1}A$ the index $k_j$ is replaced by $k_j+2,$ 
therefore the dimension of $N_A$ equals $ \sum_{k_j < 0} (|k_j|-1).$  
It follows that the index of $U_A$ is $k+N,$ where $k = \sum k_j$ is the 
total index, see also \cite{Oh}. 
\vskip .1cm
Moreover $A$ is regular if and only if 
every partial index is greater then or equal to $-1.$
\vskip .1cm
We have  $l= d-  \ (\mbox{dimension of $N_A$}) = 2m_{-} + m_0.$ 
On the other hand $(\Psi_0)(\EAC)$ has complex dimension $m_{-} + m_0 = r.$
In particular $A$ is strongly regular if and only if every partial index 
is greater then or equal to $0.$
\vskip .1cm
Similarly we can compute $r_1$ and $l_1.$ 
\vskip .1cm
Let us look at the index of $U_A(1).$
We find $dim_{\R} \LA(1) = k_{+}. $
On the other hand the dimension of the cokernel of $U_A(1)$ is 
$d-r= \sum_{k_j < 0} |k_j|,$ therefore $U_A(1)$ has index $k.$ 

\section{The case $m=1$ (hypersurface case)} \label{hyper}

\begin{lemma}
Assume $m=1,$
given $x \in \Delta.$ 
we have the following possibilities  

\[ \begin{array}{l} 
\mbox{ \ $d=l=r=0,$ \ hence A is regular and $W_x= \C^N.$}  \\ 
\mbox{$d=l=r=1,$ hence  A is regular and $W_x$ is a real hyperplane in 
$\R^{2N}.$} \\
\mbox{$d > 1, \ l=2, \ r=1$ and A is not regular.} \\
\mbox{In this case $W_x$ is a complex hyperplane in $\C^N.$} \\ 
\mbox{In particular we can not have $d=2.$} \end{array} \]
\label{m2}  \end{lemma}

\begin{proof}
If $d=0,$ then $l=r=0.$ If $d \neq 0,$ by the Proposition 
\ref{properties} $r=1$ and $l$ is either 1 or 2. 
Moreover, by 
Proposition \ref{grasm1}, and Example \ref{strongreg}, $l=1,$ if and only if  
$A$ is regular, if and only if $d=l.$ 
If we had $d=2,$ we would have $r=1.$ If it was $l=1,$ we 
would have $l=r,$ so $A$ regular, so $l=d=2.$ If it was $l=2,$ we 
would have $l=d,$ hence $A$ regular, and  $l=r=1.$ So we found an 
absurd. The other statements follow from Theorem \ref{val0}.   \end{proof}

More generally we have the following

\begin{proposition}
If $m=1,$ then $d$ is either zero or an odd positive   integer.
For every d which is either 0 or an odd positive integer, 
there exists an analytic disc attached to an hypersurface in $\C^2$
with defect d. \label{disp} \end{proposition}

\begin{proof}
Let   $2k_0$ be  the smallest strictly positive even  
integer such that there exists an $A$ with $m=1$ and $d=2k_0.$ 
From the above Lemma we know that $2k_0 > 2,$ and that $l=2,$ so 
$2k_0-2$ is even and strictly positive. However 
$2k_0-2$ is the defect of $ {\zeta^{-1}A}.$ 
This is against the minimality of $2k_0.$ 

Let us consider the example in \cite{Pa} given by the hypersurface
\[ S= \{ (z_1,z_2) \in \C^2 : Re(z_1^k z_2)=0, z_1 \neq 0 \} \]
and the disc $f(\zeta)= (\zeta,0)$ attached to  S. 
If $ k < 0$ the disc  $f$ has defect  0, while if  $k \geq 0,$  the disc has 
defect $2k+1.$ \end{proof}

Moreover

\begin{proposition}
Assume that $m=1.$
Given $y \in \de \Delta$ and $x \in  \Delta, $ we have that  
$W_x(y)$ is a  complex space, more precisely 
\vskip .1cm
If $\EA =0,$  then $W_x(y)= \C^N.$ 
\vskip .1cm
If $\EA \neq 0,$ then $W_x(y)$ is a complex 
hyperplane in $\C^N.$ \end{proposition}

\begin{proof}
By Proposition \ref{properties} we know that  $r \leq 1.$ If $r=1$ we 
conclude applying  Proposition \ref{r=m}, if $r=0,$ then $\EA=0,$ and so 
$\zeta A $ is regular. It follows from example \ref{strongreg} that 
$ \zeta A$ is strongly regular, so $l_1= r_1$ and by 
Theorem \ref{val2} we conclude. \end{proof}

\begin{proposition}
Assume that $m=1.$
Fix a point  $x \in  \de \Delta. $ 
\vskip .1cm
If $A$ is not  regular, then, $W_x= Ker(A(x)).$
\vskip .1cm
If A is regular, then $W_x= Ker(Re(A(x))).$ \end{proposition}

\begin{proof}
$A$ is regular, if and only if its strongly regular if and only if  
$l=r. $ (See Proposition \ref{grasm1}). We conclude using Theorem 
\ref{val1}. \end{proof}  

\begin{proposition}
Assume that $m=1.$
Fix two points  $x$ and $y$ in $ \de \Delta $ with $x \neq y.$   
\vskip .1cm
If $\EA \neq 0,$ then $W_x(y)= Ker(A(x))$
\vskip .1cm
If $\EA =0,$  then $W_x(y)= Ker(Re(A(x)).$ \end{proposition}

\begin{proof}
It follows directly from Theorem \ref{val3}. \end{proof}

Suppose now that $A_{f_0}$ is the map associated to an analytic disc $f_0$ 
attached to an hypersurface $S$ in $\C^N.$ Let $q=f_0(1).$ 

\begin{corollary}
If $\EA=0$  and $m=1,$  then there exists a neighborhood  $U$  
of $f_0$ in $M_q$ such that  ${\mu}_0(U)$ is an open set in 
$\C^N.$ \end{corollary}

\begin{proof}
Since $\zeta A$ is strongly regular if and only if it is regular, if and 
only if $\EA=0,$ the result follows from Theorem \ref{openmap} \end{proof}

Stefano Trapani
\vskip .1cm
Dipartimento di Matematica
\vskip .1cm
Universita' di Roma2 Tor Vergata
\vskip .1cm
Via della Ricerca Scientifica
\vskip .1cm
00133 Roma Italy
\vskip .1cm
e mail address 
\vskip .1cm
trapani@mat.utovrm.it

\end{document}